\documentclass[]{article}
\usepackage{amsmath}
\usepackage{amssymb}
\title{Infinite Product Representation for the Szeg\"{o} Kernel for an Annulus}
\author{Nuraddeen S. ~Gafai\footnote{Department of Mathematics and Statistics, Umaru Musa Yar'adua University Katsina;
		Department of Mathematical Sciences, Faculty of Science, Universiti Teknologi Malaysia, 81310 UTM Johor Bahru, Johor, Malaysia}, Ali H. M. ~Murid\footnote{Department of Mathematical Sciences, 
		Faculty of Science, Universiti Teknologi Malaysia,
		81310 UTM Johor Bahru, Johor, Malaysia}, Nur H. A. A.~Wahid\footnote{Faculty of Computer and Mathematical Sciences,
		Universiti Teknologi MARA,
		40450 Shah Alam, Selangor, Malaysia}}

\begin{document}

\maketitle

\begin{abstract}
The Szeg\"{o} kernel has many applications to problems in conformal mapping and satisfies the Kerzman-Stein integral equation. The Szeg\"{o} kernel for an annulus can be expressed as a bilateral series. In this paper, we show how to represent the Szeg\"{o} kernel for an annulus as a basic bilateral series and an infinite product representation through the application of the Ramanujan's sum. The infinite product clearly exhibits the zero of the Szeg\"{o} kernel for an annulus. Its connection with basic gamma function and modified Jacobi theta function is also presented. The results are extended to the Szeg\"{o} kernel for general annulus and weighted Szeg\"{o} kernel. Numerical comparisons on computing the Szeg\"{o} kernel for an annulus based on the bilateral series and the infinite product are also presented.
%\PACS{30D40 \and  33D15 \and  30C15 \and  65E05 \and  65R20}
\end{abstract}

\section{Introduction}
The Ahlfors map is a branching $n$-to-one map from an $n$-connected region onto the unit disk. It is intimately tied to the Szeg\"{o} kernel of an $n$-connected region \cite{S. Bell}. The boundary values of the Szeg\"{o} kernel satisfy the Kerzman-Stein integral equation, which is a Fredholm integral equation of the second kind for a region with smooth boundary \cite{BELL2}. The boundary values of the Alhfors map are completely determined from the boundary values of the Szeg\"{o} kernel \cite{S. Bell,BELL2,T}. For an annulus region $\Omega$, the Szeg\"{o} kernel can be expressed as a bilateral series from which the zero can be determined analytically \cite{TT}. The Kerzman-Stein integral equation has been solved using the Adomian decomposition method in \cite{HAZ} to give another bilateral series form for the Szeg\"{o} kernel for $\Omega$ that converges faster. There are various special functions in the form of bilateral and basic bilateral series \cite{Slater,George,Warren}. For example, the bilateral basic hypergeometric series contain, as special cases, many interesting identities related to infinite products, theta functions and Ramanujan identities. It is therefore natural to ask if the bilateral series for the Szeg\"{o} kernel for $\Omega$ can be summed as special functions or infinite product that exhibits clearly its zero.

In this paper, we show how to express the bilateral series for the Szeg\"{o} kernel for $\Omega$ as a basic bilateral series (also known as $q$-bilateral series). The Ramanujan's sum is then applied to obtain the infinite product representation for the Szeg\"{o} kernel for $\Omega$. The product clearly exhibits the zero of the Szeg\"{o} kernel for $\Omega$ and its connection with the $q$-gamma function and the modified Jacobi theta function is shown. Using symmetry of the Ramanujan's sum, we show how to easily transform the bilateral series for the Szeg\"{o} kernel for $\Omega$ in \cite{TT} to the bilateral series in \cite{HAZ}.

The plan of the paper is as follows: After the presentation of some preliminaries in Section 2, we derive the basic bilateral series and infinite product representations for the Szeg\"{o} kernel for $\Omega$ in Section 3. We then derive a closed form of the Szeg\"{o} for $\Omega$ in terms of $q$-gamma function and the modified Jacobi theta function. In Section 4, we show how to extend the representations in Section 3 to general annulus using the transformation formula for the Szeg\"{o} kernel under conformal mappings. Similar $q$-analysis for the weighted Szeg\"{o} kernel for $\Omega$ is presented in Section 5. In Section 6, we give numerical comparisons for computing the Szeg\"{o} kernel for $\Omega$ using bilateral series, infinite product, and integral equation formulations.

\section{Preliminaries}
\label{sec:1}
Let $\Omega=\{z:\rho<\mid z \mid< 1\}$ be an annulus with $0<\rho<1$ and a point $a\in\Omega$. The boundary $\Gamma$ of $\Omega$ consists of two smooth Jordan curves with outer curve $\Gamma_{0}$ oriented counterclockwise and the inner curve $\Gamma_{1}$ oriented clockwise. The positive direction of the contour $\Gamma=\Gamma_{0}\cup\Gamma_{1}$ is usually that for which the region is on the left as one traces the boundary.

Let $\{\varphi_{n}(z)\}_{n=1}^\infty$ be an orthonormal basis for the Hardy spaces $H^2(\Gamma)$. Since the Szeg\"{o} kernel $S(z,a)$ is the reproducing kernel for $H^2(\Gamma)$, it can be written as \cite{TT} 
\begin{equation}\label{1}
	S(z,a)=\sum\limits_{n=0}^\infty\varphi_n(z)\overline{\varphi_n(a)},\quad a\in \Omega,
\end{equation}
with absolute and uniform convergence on compact subsets of $\Omega$. An orthogonal basis for $H^2(\Gamma)$  is $\{z^{n}\}_{n=-\infty}^\infty$. Thus 
\begin{equation}\label{2}
	\Vert z^n\Vert^2=\int_\Gamma \vert z \vert ^{2n}\vert dz\vert = 2\pi(1+\rho^{2n+1}),
\end{equation}
where $\vert dz \vert$ is the arc length measure. Therefore an orthonormal basis for  $H^2(\Gamma)$ is \cite{T,TT}
\begin{equation}\label{3}
	\left\{\frac{z^n}{\sqrt{2\pi(1+\rho^{2n+1})}}\right\}_{n=-\infty}^\infty.
\end{equation}
Using \eqref{1} and \eqref{3} the series representation for the Szeg\"{o} kernel for $\Omega$ is given by \cite{TT}
\begin{equation}\label{4}
	S(z,a)= \frac{1}{2\pi} \sum\limits_{n=-\infty}^\infty \frac{(z\overline{a})^n}{1+\rho^{2n+1}}, \quad a\in \Omega,\quad z\in\Omega\cup\Gamma.
\end{equation}
The series \eqref{4} is a bilateral series. It has a zero at $z=-\rho/\overline{a}$ \cite{TT}.

Another bilateral series representation for the Szeg\"{o} kernel for $\Omega$ is given by \cite{HAZ} (in an equivalent form)
\begin{equation}\label{5}
	S(z,a)= \frac{1}{2\pi} \sum\limits_{n=-\infty}^\infty \frac{(-1)^{n}\rho^{n}}{\rho^{2n}-z\overline{a}}, \quad z\in\Omega\cup\Gamma, \quad\ a\in \Omega,
\end{equation}
which is initially obtained by solving the Kerzman-Stein integral equation using the Adomian decomposition method. It is also shown in \cite{HAZ} how to derive \eqref{5} directly from \eqref{4} using geometric series. It is illustrated in \cite{HAZ} that the series \eqref{5} converges faster than \eqref{4}.

More generally, if $\Omega_1$ is any doubly connected region with smooth boundary $\Gamma_1$, and $f(z)$ is a biholomorphic map of  $\Omega_1$ onto $\Omega$, then the Szeg\"{o} kernel for $\Omega_1$ can be obtained via the transformation formula as \cite{S. Bell}

\begin{eqnarray}
		S_1(z,a) &=& \sqrt{f'(z)}S(f(z), f(a))\overline{\sqrt{f'(a)}}\nonumber\\
		&=&\frac{\sqrt{f'(z)}\sqrt{f'(a)}}{2\pi}\sum\limits_{n=-\infty}^{\infty}\frac{(f(z)\overline{f(a)})^n}{1+\rho^{2n+1}},\quad a\in\Omega_1,\quad z\in\Omega_1\cup\Gamma_1,\label{6}
\end{eqnarray}
where $\rho$ is unknown but can be computed.

The Szeg\"{o} kernel $S_1(z,a)$ can also be computed without using conformal mapping. The boundary values of the Szeg\"{o} kernel $S_1(z,a)$ on $\Gamma_1$ satisfy the Kerzman-Stein integral equation \cite{{BELL2},{TT}},
\begin{equation}\label{4a}
	S_1(z,a)+\int_\Gamma A(z,w)S_1(w,a)\vert dw\vert =g(z), \quad z\in\Gamma_1 ,
\end{equation}
where 
\begin{equation*}\label{4b}
	A(z,w)=\left\{
	\begin{array}{ccc}
		\displaystyle\frac{1}{2\pi}\displaystyle\left(\frac{T(w)}{z-w}-\frac{\overline{T(z)}}{\overline{z}-\overline{w}}\right),          & z\ne w\in\Gamma_1,\\[4ex]
		0,            
		&\quad z=w\in\Gamma_1 ,
	\end{array}\right.
\end{equation*}
\begin{equation*}
	g(z)=-\frac{1}{2\pi i}\frac{\overline{T(z)}}{\overline{z}-\overline{a}}, \quad z\in\Gamma_1,
\end{equation*}
\begin{equation*}
	T(z)=\frac{z'(t)}{\vert z'(t)\vert}, \quad z\in\Gamma_1,
\end{equation*}
and $z(t)$ is a parametrization of $\Gamma_1$. The function $A(z,w)$ is known as the Kerzman-Stein kernel and it is continuous on the boundary of $\Omega_1$ \cite{Kerzman1,Kerzman2}.  In fact the integral equation \eqref{4a} is also valid for an $n$-connected region.

Since bilateral series and basic bilateral series will be used throughout this paper, we  recall some facts about $q$-series notations and results. 

Let $0<q<1$ and $\alpha\in\mathbb{C}$. The $q$-shifted factorial is defined as \cite{George}
\begin{equation}\label{7}
	(q^\alpha ;q)_n=\left\{
	\begin{array}{ccc}
		1,          & n=0,\\[4ex]
		(1-q^\alpha)(1-q^{\alpha+1})\cdots(1-q^{\alpha+n-1}),            
		&\quad n=1,2,\ldots,\\[3ex]
		\displaystyle\frac{1}{(1-q^{\alpha -1})(1-q^{\alpha -2})\cdots(1-q^{\alpha -n})}, 
		&\quad n=-1,-2,\ldots.
	\end{array}\right.
\end{equation}
This notation yields the shifted factorial as a special case through
\begin{equation}\label{8}
	\lim_{q\to1}\frac{(q^\alpha ;q)_n}{(q;q)_n}=\alpha (\alpha +1)\cdots(\alpha +n-1), \quad n=1,2,\ldots.
\end{equation}
If $\alpha$ is written in place of $q^{\alpha}$, then \eqref{7} becomes
\begin{equation}\label{9}
	(\alpha ;q)_n=\left\lbrace
	\begin{array}{ccc}
		1,          
		& n=0,\\[4ex]
		\displaystyle(1-\alpha)(1-\alpha q)\cdots(1-\alpha q^{n-1}),            
		&\quad n=1,2,\ldots,\\[3ex]
		\displaystyle\frac{1}{(1-\alpha q^{-1})(1-\alpha q^{-2})\cdots(1-\alpha q^{-n})}, 
		&\quad n=-1,-2,\ldots .
	\end{array}\right.
\end{equation}
It can be shown that \cite{George},
\begin{equation}\label{10}
	\frac{1-\alpha}{1-\alpha q^{n}}=\frac{(\alpha;q)_n}{(\alpha q;q)_n}, \quad n=0,\pm1,\pm2,\ldots.
\end{equation}
If $n\to\infty$, it is standard to write
\begin{equation}\label{11}
	(\alpha;q)_\infty=\prod_{n=0}^{\infty}(1-\alpha q^n)
\end{equation}
which is absolutely convergent for all finite values of $\alpha$, real or complex, when $\vert q\vert<1$ \cite{Slater}. This yields 
\begin{equation}\label{12}
	(\alpha;q)_n=\frac{(\alpha;q)_\infty}{(\alpha q^n;q)_\infty}.
\end{equation}
Observe that $(\alpha;q)_\infty$ would have zero as a factor if $\alpha=1$. It would be zero also if $\alpha=q^{-1},q^{-2},q^{-3},\ldots,$ but these are all outside the circle $\vert z\vert=1$ since $\vert q\vert<1$ \cite{Warren}.

The bilateral basic hypergeometric series in base $q$ with one numerator and one denominator parameters is defined by \cite{Slater,George,Warren}
\begin{equation}\label{13}
	_1\psi_1(\alpha;\beta;q;z) = \sum\limits_{n=-\infty}^\infty\frac{(\alpha;q)_n}{(\beta;q)_n}z^n.
\end{equation}
The series is convergent for $\left\vert q \right\vert < 1$ and $\left\vert \beta/\alpha \right\vert <\left\vert z \right\vert <1$.

The classical Ramanujan's $_1\psi_1$ summation is given by \cite{George,Warren}
\begin{equation}\label{15}
	_1\psi_1(\alpha;\beta;q;z)=\frac{(\alpha z;q)_\infty(q/\alpha z;q)_\infty(\beta/\alpha;q)_\infty(q;q)_\infty}{(z;q)_\infty(\beta/\alpha z;q)_\infty(q/\alpha;q)_\infty(\beta;q)_\infty},\quad \left\vert \beta/\alpha \right\vert < \left\vert z \right\vert < 1.
\end{equation}
The special case $\beta=\alpha q$ of Ramanujan's $_1\psi_1$ summation yields \cite{Warren}
\begin{equation}\label{16}
	\sum\limits_{n=-\infty}^\infty\frac{z^n}{1-\alpha q^n}=\frac{(\alpha z;q)_\infty(\frac{q}{\alpha z};q)_\infty(q;q)_\infty^2}{(z;q)_\infty(\frac{q}{z};q)_\infty(\alpha;q)_\infty(\frac{q}{\alpha};q)_\infty},
\end{equation}
also known as Cauchy's formula. Due to symmetry in $\alpha$ and $z$ on the right-hand side of \eqref{16}, it implies \cite{Warren}
\begin{equation}\label{17}
	\sum\limits_{n=-\infty}^\infty\frac{z^n}{1-\alpha q^n}=\sum\limits_{n=-\infty}^\infty\frac{\alpha^n}{1-zq^n}.    
\end{equation}
The $q$-gamma function is defined as \cite{George}
\begin{equation}\label{14}
	\Gamma_q(x) = \frac{(q;q)_\infty}{(q^x;q)_\infty}(1 - q)^{1-x}, \quad 0<q<1,\quad x=\mathbb{C}-\{0,-1,-2,\ldots\}.
\end{equation}
Another important special function that is used in this paper is the modified Jacobi theta function defined by \cite{George} 
\begin{equation}\label{17a}
	\theta(x;q)=(x;q)_\infty(q/x;q)_\infty,
\end{equation}
where $x\ne0$ and $\vert q \vert<1$. For more detailed discussion on $q$-series and historical perspectives, see for examples \cite{Slater,George,Warren} and the references therein.

\section{Szeg\"{o} Kernel for an Annulus and Basic Bilateral Series}
In this section, we express the bilateral series \eqref{4} as a basic bilateral series and derive the infinite product representation of the Szeg\"{o} kernel for $\Omega$. It is given in the following theorem.

\newtheorem{theorem}{Theorem}
\begin{theorem}
	Let $\Omega$ be the annulus $\{z: \rho<|z|<1\}$ bounded by $\Gamma$. For $a\in \Omega$, $z\in\Omega\cup\Gamma$, the Szeg\"{o} kernel for $\Omega$  can be represented by 
	\begin{eqnarray}
			S(z,a) &=& \frac{1}{2\pi(1+\rho)}         {_1\psi_1}(-\rho;-\rho^3;\rho^2;\overline{a}z)\label{18a}\\
			&=&\frac{1}{2\pi}\prod_{n=0}^{\infty}\frac{(1+\overline{a}z\rho^{2n+1})(\overline{a}z+\rho^{2n+1})(1-\rho^{2n+2})^2}{(1-\overline{a}z\rho^{2n})(\overline{a}z-\rho^{2n+2})(1+\rho^{2n+1})^2}. \label{18b}
	\end{eqnarray}
	The zero of S(z,a) in $\Omega$  is  the zero of the factor $\overline{a}z+\rho$, that is, $z=-\rho/\overline{a}$.
\end{theorem}

\noindent{\bf Proof:}
	From \eqref{4} we have
	\begin{equation}\label{19}
		S(z,a)= \frac{1}{2\pi} \sum\limits_{n=-\infty}^\infty \frac{(\overline{a}z)^n}{1+\rho^{2n+1}}=\frac{1}{2\pi} \sum\limits_{n=-\infty}^\infty \frac{(\overline{a}z)^n}{1-(-\rho)\rho^{2n}}.
	\end{equation}
	Letting $\alpha=-\rho$, and $q=\rho^2$, yields
\begin{eqnarray}
		S(z,a)&=&\frac{1}{2\pi} \sum\limits_{n=-\infty}^\infty \frac{(\overline{a}z)^n}{1-\alpha q^n}\label{19a}\\
		&=&\frac{1}{2\pi(1-\alpha)} \sum\limits_{n=-\infty}^\infty \frac{1-\alpha}{1-\alpha q^n}(\overline{a}z)^n. \label{21}
\end{eqnarray}
	
	Applying \eqref{10} and \eqref{13}, gives
\begin{eqnarray}
		S(z,a) &=&\frac{1}{2\pi(1-\alpha)}\sum\limits_{n=-\infty}^\infty\frac{(\alpha,q)_n}{(\alpha q,q)_n}(\overline{a}z)^n\nonumber\\
		&=&\frac{1}{2\pi(1-\alpha)} \quad _1\psi_1(\alpha;\alpha q;q;\overline{a}z).\label{22}
\end{eqnarray}
	
	Note that the $_1\psi_1$ series above is convergent because $\vert q \vert=\rho^2<1$ and $\vert \beta/\alpha\vert=\vert \alpha q / \alpha\vert=\vert q \vert=\rho^2<\vert \overline{a}z \vert<1$. Substituting $\alpha=-\rho$ and $q=\rho^2$ into \ref{22} gives \eqref{18a}.
	
	Applying the Ramanujan's sum \eqref{15} to \eqref{22}, gives
	\begin{equation}\label{23}
		S(z,a) =\frac{1}{2\pi(1-\alpha)}\frac{(\alpha\overline{a}z;q)_\infty(q/\alpha\overline{a}z;q)_\infty(q;q)_\infty^2}{(\overline{a}z;q)_\infty(q/ \overline{a}z;q)_\infty(q/\alpha;q)_\infty(\alpha q;q)_\infty}.
	\end{equation}
	But from \eqref{12}, with $n=1$, we have 
	$$(1-\alpha)(\alpha q;q)_\infty=(\alpha;q)_\infty.$$
	Thus \eqref{23} becomes
\begin{eqnarray}
		S(z,a)& =&\frac{1}{2\pi}\frac{(\alpha\overline{a}z;q)_\infty(q/\alpha\overline{a}z;q)_\infty(q;q)_\infty^2}{(\overline{a}z;q)_\infty(q/ \overline{a}z;q)_\infty(q/\alpha;q)_\infty(\alpha;q)_\infty}\label{23a}\\
		&=&\frac{1}{2\pi}\prod_{n=0}^{\infty}\frac{(1-\alpha\overline{a}zq^{n})(1-q^{n+1}/\alpha\overline{a}z)(1-q^{n+1})^2}{(1-\overline{a}zq^n)(1-q^{n+1}/\overline{a}z)(1-q^{n+1}/\alpha)(1-\alpha q^{n})}.\label{23b}
\end{eqnarray}
	
	Substituting $\alpha=-\rho$ and $q=\rho^2$ into \eqref{23b} gives \eqref{18b}.
	
	The infinite product \eqref{18b} would have poles if
	\[
	1-\overline{a}z\rho^{2n}=0\quad\mathrm{or}\quad \overline{a}z-\rho^{2n+2}=0 
	\]
	which implies
	\[
	z=\frac{1}{\overline{a}\rho^{2n}}\quad\mathrm{or}\quad z=\frac{\rho^{2n+2}}{\overline{a}}. 
	\]
	But
	\[
	\frac{1}{| a\rho^{2n}|}>1, \quad 
	\left|\frac{\rho^{2n+2}}{\overline{a}}\right|<\rho^{2n+1}<\rho.
	\]
	Therefore the poles are all outside $\Omega$.
	
	The infinite product \eqref{18b} would have zeros if
	\[
	1+\overline{a}z\rho^{2n+1}=0\quad\mathrm{or}\quad  \overline{a}z+\rho^{2n+1}=0
	\]
	which implies
	\[
	z=-\frac{1}{\overline{a}\rho^{2n+1}}\quad\mathrm{or}\quad z=-\frac{\rho^{2n+1}}{\overline{a}}.
	\]
	For the first case,
	\[
	\frac{1}{\vert a\rho^{2n+1}\vert}>\frac{1}{\rho^{2n+1}}>1 
	\]
	which is outside $\Omega$. For the second case, observe that
	\[
	\rho^{2n+1}<\left|\frac{\rho^{2n+1}}{\overline{a}}\right|= \frac{\rho^{2n+1}}{\vert a \vert}<\rho^{2n},
	\]
	which clearly has a zero inside $\Omega$ when $n=0$. Thus the infinite product \eqref{18b} for $S(z,a)$ has only one zero inside $\Omega$ at $z=-\rho/\overline{a}$. This completes the proof.

We note that the series representation \eqref{18a} for $S(z,a)$ is valid only for $\rho\le |z| \le 1$, while the infinite product representation \eqref{18b} for $S(z,a)$ is meaningful for all $z\in\mathbb{C}$ except for the infinitely many poles at $z=0, \rho^{-2n}/\overline{a}, \rho^{2n+2}/\overline{a}$. 

We next show that the Szeg\"{o} kernel for $\Omega$ can also be expressed in terms of the basic gamma function and modified Jacobi theta function. By applying \eqref{17a} to \eqref{23a}, and substituting $\alpha=-\rho$ and $q=\rho^2$, we have
\begin{eqnarray}
		S(z,a)&=&\frac{1}{2\pi}\frac{\theta(\alpha\overline{a}z;q)_\infty(q;q)_\infty^2}{\theta(\overline{a}z;q)_\infty(q/\alpha;q)_\infty(\alpha;q)_\infty}\nonumber\\
		&=&\frac{1}{2\pi}\frac{\theta(-\rho\overline{a}z;\rho^2)_\infty(\rho^2;\rho^2)_\infty^2}{\theta(\overline{a}z;\rho^2)_\infty(-\rho;\rho^2)_\infty^2}.\label{24b}
\end{eqnarray}

Applying \eqref{14} with $q=\rho^2$, observe that
\[
\frac{(\rho^2;\rho^2)_\infty}{(-\rho;\rho^2)_\infty}=\frac{(\rho^2;\rho^2)_\infty}{(\rho^{2x};\rho^2)_\infty}=\frac{\Gamma_{\rho^2}(x)}{(1-\rho^2)^{1-x}},
\]
where $x$ satisfies $\rho^{2x}=-\rho$. This equation may be written as
\[
e^{(2x-1)\ln\rho}=e^{i\pi},
\]
which yields a solution
\[
x=\frac{1}{2}+\frac{i\pi}{2\ln\rho}.
\]
Thus \eqref{24b} becomes
\begin{equation}\label{24c}
	S(z,a)=\frac{[\Gamma_{\rho^{2}}(\lambda)]^2}{2\pi (1-\rho^2)^{2(1-\lambda)}}\frac{\theta(-\rho\overline{a}z;\rho^2)_\infty}{\theta(\overline{a}z;\rho^2)_\infty}, \quad\lambda=\frac{1}{2}+\frac{i\pi}{2\ln\rho}.
\end{equation}
This can be regarded as a closed form expression for the Szeg\"{o} kernel for $\Omega$. 

In the following we show how to easily transform the series \eqref{4} to the series \eqref{5} using \eqref{17}. Letting $\alpha=-\rho$ and $q=\rho^2$, \eqref{4} becomes

\begin{equation*}
	S(z,a)=\frac{1}{2\pi }\sum\limits_{n=-\infty}^\infty\frac{(\overline{a}z)^{n}}{1-\alpha q^n} =\frac{1}{2\pi }\sum\limits_{n=-\infty}^\infty\frac{\alpha^n}{1-(\overline{a}z)q^n},		
\end{equation*}
where in the last step we have used \eqref{17}. By replacing $\alpha=-\rho$ and $q=\rho^2$, we get

\begin{equation*}\label{30}
	S(z,a)=\frac{1}{2\pi }\sum\limits_{n=-\infty}^\infty\frac{(-1)^{n}\rho^n}{1-(\overline{a}z)\rho^{2n}}.	
\end{equation*}
Letting $n=-m$, yields

\begin{equation*}
	S(z,a)=\frac{1}{2\pi }\sum\limits_{m=-\infty}^\infty\frac{(-1)^{-m}\rho^{-m}}{1-(\overline{a}z)\rho^{-2m}}
	=\frac{1}{2\pi }\sum\limits_{m=-\infty}^\infty\frac{(-1)^{m}\rho^m}{\rho^{2m}-\overline{a}z},
\end{equation*}
which is the same as \eqref{5}.

\section{Szeg\"{o} Kernel for General Annulus}
Consider the general annulus $\Omega_2=\left\{z:r_2<\mid z-z_0 \mid<r_1\right\}$ with boundary denoted by $\Gamma_2$. The region $\Omega_2$ reduces to $\Omega$ if $z_0=0$, $r_2=\rho$ and $r_1=1.$

\begin{theorem}
	Let $z_0\in\mathbb{C}$, $z\in\Omega_2\cup\Gamma_2$ and $a\in\Omega_2$. The Szeg\"{o} kernel for $\Omega_2$ can be represented by the bilateral series as
	\begin{eqnarray}
			S_2(z,a) &=&\frac{1}{2\pi}\sum\limits_{n=-\infty}^\infty\frac{(\overline{a}-\overline{z_0})^n}{r_1^{2n+1}+r_2^{2n+1}} (z-z_0)^n\label{25}\\
			&=&\frac{1}{2\pi}\sum\limits_{n=-\infty}^\infty\frac{(-1)^nr_1^{n+1}r_2^n}{r_2^{2n}r_1^2-r_1^{2n}(z-z_0)(\overline{a}-\overline{z_0})}. \label{25b}
	\end{eqnarray}
	
	The zero of $S_2(z,a)$ in $\Omega_2$ is $z=z_0-\frac{r_1r_2}{\overline{a}-\overline{z_0}}$.
\end{theorem}
	 
\noindent{\bf Proof:}
	Observe that the function $f(z)=(z-z_0)/r_1$ maps $\Omega_2$ onto $\Omega$ with $\rho=r_2/r_1$.\\Applying the transformation formula \eqref{6}, yields
\begin{eqnarray}
		S_2(z,a) &=&\sqrt{f'(z)}S(f(z),f(a))\overline{\sqrt{f'(a)}}\nonumber\\
		&=&\frac{1}{\sqrt{r_1}}S\left(\displaystyle\frac{z-z_0}{r_1},\frac{a-z_0}{r_1}\right)\frac{1}{\overline{\sqrt{r_1}}}\nonumber\\
		&=&\frac{1}{r_1}S\left(\displaystyle\frac{z-z_0}{r_1},\frac{a-z_0}{r_1}\right).\label{26} 
\end{eqnarray}
	
	Applying \eqref{4} to \eqref{26} with $z$ and $a$ replaced by $(z-z_0)/r_1$ and $(a-z_0)/r_1$ respectively, gives
	\begin{equation}\label{27}
		S_2(z,a)=\frac{1}{2\pi r_1}\sum\limits_{n=-\infty}^\infty\frac{((z-z_0)(\overline{a-z_0})/r_1^2)^n}{1+(r_2/r_1)^{2n+1}},
	\end{equation}
	which simplifies to \eqref{25}. 
	
	Applying \eqref{5} to \eqref{26} instead with $z$ and $a$ replaced by $(z-z_0)/r_1$ and $(a-z_0)/r_1$ respectively, gives
	\begin{equation}\label{27a}
		S_2(z,a)=\frac{1}{2\pi r_1}\sum\limits_{n=-\infty}^\infty\frac{(-1)^n(r_2/r_1)^n}{(r_2/r_1)^{2n}-(z-z_0)(\overline{a-z_0})/r_1^2},
	\end{equation}
	which simplifies to \eqref{25b}.
	
	Using the fact that $S(z,a)$ has a zero at $z=-\rho/\overline{a}$ for $\Omega$, the zero of $S_2(z,a)$ for $\Omega_2$ is $\frac{z-z_0}{r_1}=\frac{-\rho}{(\overline{a-z_0})/r_1}$ which implies  $z=z_0-\frac{\rho r_1^{2}}{(\overline{a}-\overline{z_0})}=z_0-\frac{r_1r_2}{(\overline{a}-\overline{z_0})}$.
	This completes the proof.

Similarly, the infinite product representation of $S_2(z,a)$ for $\Omega_2$ can be obtained by applying \eqref{18b} to \eqref{26} with $z$ and $a$ replaced by $(z-z_0)/r_1$ and $(a-z_0)/r_1$ respectively.

\section{The Weighted Szeg\"{o} Kernel for an Annulus and Basic Bilateral Series}
The weighted Szeg\"{o} kernel is defined in \cite{Scott} as 

\begin{equation}\label{36}
	\hat{K}_{q}^{t}(z,w)=\frac{1}{2\pi}\sum\limits_{n=-\infty}^{\infty}\frac{(\overline{w}z)^n}{1+tq^{2n}}, \quad t>0, \quad q<\vert z \vert, \vert w \vert<1. 
\end{equation}
To adopt the notations used earlier, we change $q$ to $\rho$, and $w$ to $a$ and $\hat{K}_{q}^{t}(z,w)$ to $S_\rho^t(z,a)$ in \eqref{36}, which gives
\begin{equation*}
	S_{\rho}^{t}(z,a)=\frac{1}{2\pi}\sum\limits_{n=-\infty}^{\infty}\frac{(\overline{a}z)^n}{1+t\rho^{2n}}, \quad t>0, \quad \rho<\vert z \vert, \vert a \vert<1. 
\end{equation*}
Note that $S_{\rho}^{t}(z,a)$ is exactly the kernel $S(z,a)$ for $\Omega$ discussed in section 1. The zeros of the kernel $S_{\rho}^{t}(z,a)$ is not discussed in \cite{Scott} but has expressed interest on the effect of the weight on the location of its zeros. In the following theorem, we express the weighted Szeg\"{o} kernel $S_{\rho}^{t}(z,a)$ as a basic bilateral series and derive its associated infinite product representation as well as its zeros.

\begin{theorem}\label{Theorem 3}
	Let $\Omega$ be the annulus $\{z: \rho<|z|<1\}$ bounded by $\Gamma$. For $a\in \Omega$, $z\in\Omega\cup\Gamma$, $t>0$, the weighted Szeg\"{o} kernel $S_{\rho}^{t}(z,a)$ for $\Omega$  can be represented by
\begin{eqnarray}
		S^{\rho}_{t}(z,a)& =& \frac{1}{2\pi(1+t)}         {_1\psi_1}(-t;-t\rho^2;\rho^2;\overline{a}z)\label{37a}\\
		&=&\frac{1}{2\pi}\prod_{n=0}^{\infty}\frac{(1+t\overline{a}z\rho^{2n})(\overline{a}z+\rho^{2n+2}/t)(1-\rho^{2n+2})^2}{(1-\overline{a}z\rho^{2n})(\overline{a}z-\rho^{2n+2})(1+\rho^{2n+2}/t)(1+t\rho^{2n})}. \label{37b}
\end{eqnarray} 
	
	The kernel $S_{\rho}^{t}(z,a)$ has a zero in $\Omega$  only if $t$ takes the form $t=\rho^{\pm(2m+1)}$, $m=0, 1, 2, \dots$. In both cases, the zero is $z=-\rho/\overline{a}$.
\end{theorem}

\noindent{\bf Proof:}
	Observe that
	\begin{equation*}
		S_{\rho}^{t}(z,a)=\frac{1}{2\pi}\sum\limits_{n=-\infty}^{\infty}\frac{(\overline{a}z)^n}{1-(-t)\rho^{2n}}.
	\end{equation*}
	Letting $\alpha=-t$ and $q=\rho^2$, the above equation becomes
\begin{eqnarray*}
		S_{\rho}^{t}(z,a) =\frac{1}{2\pi}\sum\limits_{n=-\infty}^{\infty}\frac{(\overline{a}z)^n}{1-\alpha q^{n}}
\end{eqnarray*}
	
	which is exactly the same form as \eqref{19a}. Applying the result 
	\eqref{22} with $\alpha=-t$, the above equation becomes
\begin{eqnarray}
		S_{\rho}^{t}(z,a)=\frac{1}{2\pi(1+t)}  {_1\psi_1}(-t;-tq;q;\overline{a}z).\label{37}
\end{eqnarray}
	
	The series \eqref{37} is convergent because $\vert q \vert=\rho^2<1$ and $\vert \beta/\alpha \vert=\vert -tq/(-t)\vert = \vert q \vert <\rho^2<\vert \overline{a}z\vert<1$. Substituting $q=\rho^2$ gives \eqref{37a}.
	
	Applying  the result \eqref{23a} with $\alpha=-t$ to \eqref{37}, yields
	\begin{equation}\label{38}
		S^{t}_{\rho}(z,a)=\frac{1}{2\pi }\frac{(-t\overline{a}z;q)_\infty(q/(-t)\overline{a}z;q)_\infty(q;q)_\infty^{2}}{(\overline{a}z;q)_\infty(q/\overline{a}z;q)_\infty(q/(-t);q)_\infty(-t;q)_\infty}.
	\end{equation}
	Replacing $q=\rho^2$ and applying \eqref{11} gives \eqref{37b}.
	
	In the proof of Theorem 1, we have shown that the factors $(1-\overline{a}z\rho^{2n})(\overline{a}z-\rho^{2n+2})$ have no zeros in $\Omega$. The factors $(1+\rho^{2n+2}/t)(1+t\rho^{2n})$ would have zeros if
	\[
	\rho^{2n+2}/t=-1 \quad\mathrm{or}\quad t\rho^{2n}=-1.
	\]
	Since $t>0$, we conclude that the kernel $S^{t}_{\rho}(z,a)$ has no poles in $\Omega$ for any $t>0$. The factors $(1+t\overline{a}z\rho^{2n})(\overline{a}z+\rho^{2n+2}/t)$ would have zeros if
	\[
	1+t\overline{a}z\rho^{2n}=0\quad\mathrm{or}\quad  \overline{a}z+\rho^{2n+2}/t=0
	\]
	which implies
	\[
	z=-\frac{1}{t\overline{a}\rho^{2n}}\quad\mathrm{or}\quad z=-\frac{\rho^{2n+2}}{t\overline{a}}.
	\]
	
	For the first case, observe that
	\[
	\frac{1}{t\rho^{2n}}<\frac{1}{\vert t\overline{a}\rho^{2n}\vert}<\frac{1}{t\rho^{2n+1}}. 
	\]
	To have a zero in $\Omega$, we must have the condition
	\[
	\rho\le\frac{1}{t\rho^{2n}}<\frac{1}{\vert t\overline{a}\rho^{2n}\vert}<\frac{1}{t\rho^{2n+1}}\le 1 
	\]
	which means
	\[
	t\le\frac{1}{\rho^{2n+1}}\quad\mathrm{and}\quad t\ge\frac{1}{\rho^{2n+1}}.
	\]
	Hence we must have $t=\rho^{-(2n+1)}$. In this case, the zero of $S^{t}_{\rho}(z,a)$ in $\Omega$ is $z=-\rho/\overline{a}$.
	
	For the second case, observe that
	\[
	\frac{\rho^{2n+2}}{t}<\frac{\rho^{2n+2}}{\vert t\overline{a}\vert}<\frac{\rho^{2n+1}}{t}. 
	\]
	To have a zero in $\Omega$, we must have the condition
	\[
	\rho\le\frac{\rho^{2n+2}}{t}<\frac{\rho^{2n+2}}{\vert t\overline{a}\vert}<\frac{\rho^{2n+1}}{t}\le 1 
	\]
	which means
	\[
	t\le\rho^{2n+1}\quad\mathrm{and}\quad t\ge\rho^{2n+1}.
	\]
	Hence we must have $t=\rho^{2n+1}$. In this case, the zero of $S^{t}_{\rho}(z,a)$ in $\Omega$ is also $z=-\rho/\overline{a}$. This completes the proof.

The weighted Szeg\"{o} kernel can also be expressed in terms of basic gamma function and the modified Jacobi theta function. By applying \eqref{17a} to \eqref{38} with $q=\rho^2$, we have
\begin{equation}\label{40}
	S^{t}_{\rho}(z,a)=
	\frac{1}{2\pi}\frac{\theta(-t\overline{a}z;\rho^2)_\infty(\rho^2;\rho^2)_\infty^{2}}{\theta(\overline{a}z;\rho^2)_\infty(\rho^2/(-t);\rho^2)_\infty(-t;\rho^2)_\infty}.
\end{equation}
Observe that
\[
\frac{(\rho^2;\rho^2)_\infty}{(-t;\rho^2)_\infty}=\frac{(\rho^2;\rho^2)_\infty}{(\rho^{2x};\rho^2)_\infty}=\frac{\Gamma_{\rho^2}(x)}{(1-\rho^2)^{1-x}},
\]
where $x$ satisfies $\rho^{2x}=-t$. This equation may be written as
\[
2x\ln\rho=\ln(-t)=\ln|-t|+i\arg(-t)=\ln t+i\pi
\]
which yields a solution
\[
x=\frac{\ln t+i\pi}{2\ln\rho}.
\]
Observe also that
\[
\frac{(\rho^2;\rho^2)_\infty}{(-\rho^2/t;\rho^2)_\infty}=\frac{(\rho^2;\rho^2)_\infty}{(\rho^{2y};\rho^2)_\infty}=\frac{\Gamma_{\rho^2}(y)}{(1-\rho^2)^{1-y}},
\]
where $y$ satisfies $\rho^{2y}=-\rho^2/t$. This equation may be written as
\[
(2y-2)\ln\rho=\ln(-1/t)=\ln|-1/t|+i\arg(-1/t)=-\ln t+i\pi
\]
which yields a solution
\[
y=1+\frac{-\ln t+i\pi}{2\ln\rho}.
\]
Thus \eqref{40} becomes
\begin{equation}\label{41}
	S_\rho^t(z,a)=\frac{\Gamma_{\rho^{2}}(\mu)\Gamma_{\rho^{2}}(\nu)\theta(-t\overline{a}z;\rho^2)_\infty}{2\pi (1-\rho^2)^{2-\mu-\nu}\theta(\overline{a}z;\rho^2)_\infty}, \quad\mu=\frac{\ln t+i\pi}{2\ln\rho},\quad \nu=1+\frac{-\ln t+i\pi}{2\ln\rho}.
\end{equation}
This can be regarded as a closed form expression for the weighted Szeg\"{o} kernel for an annulus $\Omega$. Observe that \eqref{41} reduces to \eqref{24c} when $t=\rho$.

\section{Numerical Computation of the Szeg\"{o} Kernel for an Annulus}
In this section we compare the speed of convergence of the three formulas for computing the Szeg\"{o} kernel for $\Omega$ based on the two bilateral series \eqref{4}, \eqref{5} and the infinite product \eqref{18b}.

We consider an annulus $\Omega$ bounded by
$$\Gamma_{0}:z_0(t)=e^{it},$$ $$\Gamma_{1}:z_1(t)=\rho e^{-it},$$\\
with $0\le t \le2\pi$, $a=0.7i$ and $\rho=0.5$.

To approximate \eqref{4} numerically, we calculate
\begin{equation*}
	S(z,a)\approx S_{10}(z,a)=\frac{1}{2\pi}\sum\limits_{k=-10}^{10}\frac{(z\overline{a})^k}{1+\rho^{2k+1}},
\end{equation*}
$S_{50}$ and $S_{100}$.

To approximate \eqref{5}  numerically, we calculate
\begin{equation*}
	S(z,a)\approx S_{10}^*(z,a)=\frac{1}{2\pi}\sum\limits_{k=-10}^{10}\frac{(-1)^{k}\rho^{k}}{\rho^{2k}-z\overline{a}}
\end{equation*}
and $S_{50}^*$.

To approximate \eqref{18b}  numerically, we compute
\begin{equation*}
	S(z,a)\approx S_{15}^{**}(z,a)=\frac{1}{2\pi}\prod_{k=0}^{15}\frac{(1+\overline{a}z\rho^{2k+1})(z\overline{a}+\rho^{2k+1})(1-\rho^{2k+2})^2}{(1-z\overline{a}\rho^{2k})(z\overline{a}-\rho^{2k+2})(1+\rho^{2k+1})^2},
\end{equation*}
$S_{20}^{**}$ and $S_{25}^{**}$.

The approximations are then compared with the numerical solution of the Kerzman-Stein equation \eqref{4a}. To solve \eqref{4a} we used the Nystr\"{o}m method \cite{HAZ} with trapezoidal rule with $n$ selected nodes on each boundary component $\Gamma_0$ and $\Gamma_1$. The approximate solution is represented by $\tilde{S}_n$ where $n$ is the number of nodes. All the computations were done using MATHEMATICA 12.3. The results for the error norms are presented in Tables 1-3.
	
	\begin{table}[htbp]
		\caption{Error norms between $S_{10}$ and $\tilde{S}_n$, $S_{50}$ and $\tilde{S}_n$, $S_{100}$ and $\tilde{S}_n$}\label{4con1}
		\centering
		\begin{tabular}{lllll}
			\hline\noalign{\smallskip}
			{$\ n$} & {$\vert\vert S_{10}-\tilde{S}_n\vert\vert_\infty$} & {$\vert\vert S_{50}-\tilde{S}_n\vert\vert_\infty$} & {$\vert\vert S_{100}-\tilde{S}_n\vert\vert_\infty$} \\
			\noalign{\smallskip}\hline\noalign{\smallskip}
			{$16$}  & {$2.4536(-02)$}  & {$2.97754(-03)$}  & {$2.97758(-03)$}\\
			
			{$32$}  & {$2.75019(-02)$} & {$1.15906(-05)$}  & {$1.16299(-05)$}\\
			
			{$64$}  & {$2.75136(-02)$} & {$3.91113(-08)$}  & {$1.88349(-10)$}\\
			
			{$128$}  & {$2.75136(-02)$} & {$3.92996(-08)$}  & {$2.28878(-15)$}\\
			\noalign{\smallskip}\hline
		\end{tabular}
	\end{table}
	\begin{table}[htbp]
		\caption{Error norms between $S_{10}^*$ and $\tilde{S}_n$, $S_{50}^*$ and $\tilde{S}_n$}\label{4con2}
		\centering
		\begin{tabular}{lllll}
			\hline\noalign{\smallskip}
			{$\ n$} & {$\vert\vert S_{10}^*-\tilde{S}_n\vert\vert_\infty$} & {$\vert\vert S_{50}^*-\tilde{S}_n\vert\vert_\infty$} \\
			\noalign{\smallskip}\hline\noalign{\smallskip}
			{$16$}  & {$2.94797(-03)$}  & {$2.97758(-03)$}\\
			
			{$32$}  & {$1.78995(-02)$} & {$1.16299(-05)$}\\
			
			{$64$}  & {$1.77628(-04)$} & {$1.88351(-10)$}\\
			
			{$128$}  & {$1.77628(-04)$} & {$1.81497(-15)$}\\
			\noalign{\smallskip}\hline
		\end{tabular}
	\end{table}
	\begin{table}[htbp]
		\caption{Error norms between $S_{15}^{**}$ and $\tilde{S}_n$, $S_{20}^{**}$ and $\tilde{S}_n$, $S_{25}^{**}$ and $\tilde{S_n}$}\label{4con3}
		\centering
		\begin{tabular}{lllll}
			\hline\noalign{\smallskip}
			{$\ n$} & {$\vert\vert S_{15}^{**}-\tilde{S}_n\vert\vert_\infty$} & {$\vert\vert S_{20}^{**}-\tilde{S}_n\vert\vert_\infty$} & {$\vert\vert S_{25}^{**}-\tilde{S}_n\vert\vert_\infty$} \\
			\noalign{\smallskip}\hline\noalign{\smallskip}
			{$16$}  & {$2.97758(-03)$}  & {$2.97758(-03)$}  & {$2.97758(-03)$}\\
			
			{$32$}  & {$1.16296(-05)$} & {$1.16299(-05)$}  & {$1.16299(-05)$}\\
			
			{$64$}  & {$1.44308(-10)$} & {$1.88038(-10)$}  & {$1.8835(-10)$}\\
			
			{$128$}  & {$3.1999(-10)$} & {$3.1275(-13)$}  & {$1.82618(-15)$}\\
			\noalign{\smallskip}\hline
		\end{tabular}
	\end{table}
\newpage
The numerical results presented in Tables 1-3 show that computations using the infinite product formula \eqref{18b} converges faster than the bilateral series formulas \eqref{4} and \eqref{5}.

\section{Conclusion}
This paper has shown that the bilateral series for the Szeg\"{o} kernel for annulus $\Omega$ is a disguised bilateral basic hypergeometric series $_1\psi_1$. The Ramanujan's sum for $_1\psi_1$ is then applied to obtain the infinite product representation for the Szeg\"{o} kernel for annulus $\Omega$. The product clearly exhibits the zero of the Szeg\"{o} kernel for an $\Omega$. The Szeg\"{o} kernel can also be expressed as a closed form in terms of the $q$-gamma function and the modified Jacobi theta function. Similar $q$-analysis has also been conducted for the Szeg\'{o} kernel for general annulus $\Omega_2$, and for the weighted Szeg\"{o} kernel for $\Omega$. The numerical comparisons have shown that the infinite product method converges faster than the bilateral series methods for computing the Szeg\"{o} kernel for $\Omega$.

\section*{Acknowledgements}
This study was supported partially by the Ministry of Education Malaysia (MOE) through the Research Management Centre (RMC), Universiti Teknologi Malaysia (FRGS Ref. No. R.J130000.7854.5F198). The first author would also like to acknowledge the Tertiary Education Trust Fund (tetfund) Nigeria for overseas scholarship award. 

\bibliographystyle{unsrt}  


%
\end{document}